\documentclass[10pt,reqno]{amsart}
     \makeatletter
     \def\section{\@startsection{section}{1}%
     \z@{.7\linespacing\@plus\linespacing}{.5\linespacing}%
     {\bfseries
     \centering
     }}
     \def\@secnumfont{\bfseries}
     \makeatother
\newtheorem{thm}{Theorem}[section]

\newtheorem{lem}[thm]{Lemma}

\theoremstyle{definition}
\newtheorem{defn}[thm]{Definition}
\theoremstyle{remark}
\newtheorem{rem}[thm]{Remark}
\numberwithin{equation}{section} 

\newcommand{\Div}{\mathop{\mathrm{Div}}}
\renewcommand{\d}{\/\mathrm{d}\/}

\begin{document}

\title[Lyapunov Functionals and Local Dissipativity] {Lyapunov Functionals and Local
Dissipativity for the Vorticity Equation in $\mathrm{L^{p}}$ and
Besov Spaces}

\author[Utpal Manna]{Utpal Manna*\\ Department of Mathematics, University of Wyoming,\\ Laramie, Wyoming 82071,
USA\\ e-mail: utpal@uwyo.edu\\ \\}
\thanks{* This research is supported by Army Research Office, Probability and
Statistics Program, grant number DODARMY1736}

\author[S.S. Sritharan]{\\ \\S.S. Sritharan*\\Department of Mathematics, University of Wyoming,\\ Laramie, Wyoming
82071, USA\\ e-mail: sri@uwyo.edu}



\begin{abstract}
In this paper we establish the local Lyapunov property of certain
$\mathrm{L}^{p}$ and Besov norm of the vorticity fields. We have
resolved in part, certain open problem posed by Tosio Kato for the
three dimensional Navier-Stokes equation by studying the vorticity
equation. The local dissipativity of the sum of linear and
non-linear operators of the vorticity equation is established. One
of the main techniques used here is the Littlewood-Paley analysis.
\end{abstract}

\maketitle


\section{Introduction}

Stability and control of a dynamical system is often studied using
Lyapunov functions~\cite{FrKo96, La66, La62}. The local Lyapunov
property we study in this paper can thus be of interest to the
understanding, control and stabilization of turbulent
fields~\cite{Sr98}. This property also sheds some light towards the
research on global Navier-Stokes solutions in super-critical spaces
(for definitions and examples of these spaces see~\cite{Ca95}
and~\cite{CaMe95}). Weak solutions of the Navier-Stokes equation
satisfy the energy inequality which in turn implies that the
$\mathrm{L}^2$-norm of velocity decreases in time~\cite{La69}. This
idea was generalized by Tosio Kato~\cite{Ka90} to prove that for
every solution of Navier-Stokes equation in $\mathbb{R}^m$ ($m\geq
3$), there exist a large number of Lyapunov functions, which
decrease monotonically in time if the solution have small
$\mathrm{L}^m(\mathbb{R}^m)$-norm. More specifically Kato proved
that the local Lyapunov property in $\mathrm{L}^p$-norm for
$1<p<\infty$ and in $\mathrm{W}^{s,p}$-norm for $ s>0,\ 2\leq
p<\infty$. He also noted that for any Lyapunov function
$\mathfrak{L}(u)$ and any monotone increasing function $\Phi$,
$\Phi(\mathfrak{L}(u))$ will again be a Lyapunov function. Moreover
Kato also proved the local dissipativity of the sum of the linear
and nonlinear operators of the Navier-Stokes equation in
$\mathrm{L}^p$-norm for $2\leq p<\infty$. However the local
dissipativity in $\mathrm{W}^{s,p}$-norm for $s
> 0$ has remained an open problem.

Cannone and Planchon~\cite{CaPl00} proved the Lyapunov property for
the $3$-D Navier-Stokes equation in Besov spaces. In particular they
proved that if $p, q \geq 2$, $\frac{3}{p}+\frac{2}{q}>1$ and as
long as the $\dot{B}_{\infty}^{-1, \infty}$-norm of the velocity is
small, the $\dot{B}_{p}^{-1+3/p, q}$-norm of velocity decreases in
time.

In~\cite{KoTa01} Koch and Tataru considered the local and global (in
time) well-posedness for the incompressible Navier-Stokes equation
and proved the existence and uniqueness of global mild solution in
$BMO^{-1}$ provided that the initial solution is small enough in
this space. Due to Cannone and Planchon~\cite{CaPl00}, existence of
the Lyapunov functions for small $\dot{B}_{\infty}^{-1,
\infty}$-norm is known but the global solvability of Navier-Stokes
equation in this space remains an open problem. Noting the embedding
theorem $BMO^{-1}\subset \dot{B}_{\infty}^{-1, \infty}$, $BMO^{-1}$
is thus the largest space of initial data for which global mild
solution has been shown to exist.

Recently, P. G. Lemari\'{e}-Rieusset~\cite{Lr02} has extended the
result of Cannone and Planchon~\cite{CaPl00} to a larger class of
initial data. He proved that for initial data $u_{0} \in
\dot{B}_{p}^{s,q} \cap BMO^{-1}$ where $s
> -1$, $p\geq 2$, $q\geq 1$ and $s + \frac{2}{q} > 0$, there exists
a constant $C_{0}
> 0$ independent of $p$ and
$q$, such that if $u$ is a Koch-Tataru solution of Navier-Stokes
equation and satisfying $sup_{t}\parallel u(t)
\parallel_{\dot{B}_{\infty}^{-1, \infty}} < C_{0}$, then
$ t\rightarrow\parallel u(t)\parallel_{\dot{B}_{p}^{s,q}}$ is a
Lyapunov function.

Local monotonicity of different type has been used in proving the
solvability in unbounded domains for Navier-Stokes in
$2$-D~\cite{MeSr00}, in $3$-D~\cite{BaSr01} and for modified $2$-D
Navier-Stokes with artificial compressibility~\cite{MaMeSr06}. Local
monotonicity has also been useful in Control theory~\cite{BaSr01}.
For extensive theories and applications on dissipative and accretive
operators see Barbu~\cite{Ba76} and Browder~\cite{Br76}.

In this paper, we achieve a partial resolution to the open problems
posed by Kato~\cite{Ka90} for the Navier-Stokes equation by studying
the vorticity equation:
\begin{equation}\label{e1.1}\left\{\begin{aligned}
&\partial_{t} \omega\ - \nu\triangle\omega+ u\cdot\nabla\omega -
\omega\cdot\nabla u = 0,\
\text{in} \ R^{m}\times R_{+},\\
&\nabla\cdot\omega =0,\ \text{in} \ R^{m}\times R_{+},\\
&\omega(x,0)= \omega_{0}(x), \ x \in \ R^{m}.
\end{aligned}\right.\end{equation}

To be specific, we have proved that the vorticity equation have a
family of Lyapunov functions in $\mathrm{L}^{p}(\mathbb{R}^m)$ for
$2\leq p<\infty$ and $m \geq 3$ provided that the
$\mathrm{L}^{m}$-norm of the velocity is small enough. We then prove
$\dot{B}_{p}^{-1+3/p, q}$-norm of the vorticity is a Lyapunov
function for $3$-D vorticity equation provided the velocity and the
vorticity are small in $\dot{B}_{\infty}^{-1, \infty}$-norm and
$\dot{B}_{\infty}^{-2, \infty}$-norm respectively.

We have also proved the dissipativity of the sum of the linear and
nonlinear operators of the vorticity equation \eqref{e1.1} in
$\mathrm{L}^{p}$ for $2\leq p<\infty$, which in part answers the
open problem of Kato for the local dissipativity of the
Navier-Stokes operators in $\mathrm{W}^{1,p}$-norm.

In Section 2 and 3 we recall some basic facts concerning
Littlewood-Paley decomposition, homogeneous Besov spaces and the
Paraproduct rule. The main results are presented in section 4.


\section{Some Definitions and Estimates}
\begin{defn}($Duality\ Map$)
 The mapping $G: X\rightarrow 2^{X^{\star}}$ is
 called the duality mapping of the space X if $$ G(x) =
\{x^{\star}\in X^{\star};\ \langle x, x^{\star}\rangle =\ \parallel
x \parallel_{X}^{2}\ =\ \parallel x^{\star}
 \parallel_{X^{\star}}^{2},\ \forall x\in X \}.$$
\end{defn}
\begin{rem}
The duality map for $\mathrm{L}^{p}$ is given by
$$G(x)= \frac{x\mid
x\mid^{p-2}}{\parallel x\parallel_{p}^{p-2}}.$$
\end{rem}

\begin{defn}($Dissipative\ Operator$)
An operator $\textit{A}$ is said to be dissipative
 if $$\langle \textit{A} x - \textit{A} y, \ G(x-y) \rangle \leq 0,\qquad \ \forall x, y \in \emph{D}(\textit{A}).$$
 An operator $\textit{A}$ is said to be accretive if $- \textit{A}$ is
 dissipative.\\
 See~\cite{Ba76} and~\cite{Br76} for extensive theories and applications on nonlinear
 operators in Banach spaces.
\end{defn}

\begin{defn}($Lyapunov\ Function$)
Let $v$ be a solution of the Navier-Stokes equation. Then any
function $\mathfrak{L}(v)(t)$ montonically decreasing in time is
called a Lyapunov function associated to $v$.
\end{defn}
The most well-known example is certainly provided by
 energy~\cite{La69}
 $$ E(v)(t) = \frac{1}{2}\parallel v(t)\parallel_{2}^{2}.$$
The energy equality for the Navier-Stokes equation yield
\begin{align}
\frac{d}{dt}E(t) + \nu\parallel\nabla v(t)\parallel_{2}^{2}\ =
0,\nonumber
\end{align}
which proves that $E(t)$ is Lyapunov functional.

Let us now recall two lemmas due to Kato~\cite{Ka90}.

\begin{lem}
Let $2\leq p <\infty\ \text{and}\ \phi\in \mathrm{W}^{1,p}$. Define
\begin{align}\label{e2.1}
&Q_{p}(\phi) = \int_{\partial\phi(x)\neq
0}\mid\phi(x)\mid^{p-2}\mid\nabla\phi(x)\mid^{2}\d x\ \geq 0.
\end{align}
Then
\begin{align}\label{e2.2}
 CQ_{p}(\phi) \leq -\langle\mid\phi\mid^{p-2}\phi,
\Delta\phi\rangle \ < \infty,
\end{align}
where $C$ denotes a positive constant.
\end{lem}

\begin{lem}
Let $2\leq p <\infty$ and $\phi\in \mathrm{W}^{1,p}$. Then
\begin{align}\label{e2.3}
 \parallel \phi
\parallel_{\frac{mp}{m-2}} \ \leq CQ_{p}(\phi)^\frac{1}{p}.
\end{align}
\end{lem}

\begin{lem}
Let $u$ be the velocity field obtained from $\omega$ via the
Biot-Savart law:
\begin{align}\label{e2.4}
u(x) =
-\frac{\Gamma(m/2+1)}{m(m-2)\pi^{m/2}}\int_{\mathbb{R}^{m}}\frac{(x-y)}{\mid
x-y
 \mid^{m}}\ \times \omega(y) \d y, \qquad x \in {\mathbb{R}^{m}},
 m\geq 3.
\end{align}
\textbf{(a)} Assume that $ 1 < p <\infty $. Then for every
divergence-free vector field $u$ whose gradient is in
$\mathrm{L}^{p}$, there exists a $ C > 0 ,$ depending on $p$, such
that
\begin{align}\label{e2.5}
 \parallel\nabla u \parallel_{p}\leq \ \
C\parallel\omega\parallel_{p}.
\end{align}
\textbf{(b)} If $\omega\in \mathrm{L}^{1}(\mathbb{R}^{m}) \cap
\mathrm{L}^{p}(\mathbb{R}^{m})$, $\frac{m}{m-1} \ < \ p \ \leq \
\infty $, then
\begin{align}\label{e2.6}
\parallel u \parallel_{\mathrm{L}^{p}(\mathbb{R}^{m})} \ \leq
C\big(\parallel\omega\parallel_{\mathrm{L}^{1}(\mathbb{R}^{m})} \ +
\
\parallel\omega\parallel_{\mathrm{L}^{p}(\mathbb{R}^{m})}\big).
\end{align}

\end{lem}
\begin{proof}
\textbf{(a)} See Theorem 3.1.1 in~\cite{Ch98}.\\
\textbf{(b)} The proof is due to Ying and Zhang~\cite{LyZp96}, Lemma
3.3.1.
\end{proof}


\section{Littlewood-Paley Decomposition and \\ Besov Spaces}

In this section, we recall some classical results concerning the
homogeneous Besov spaces in terms of the Littlewood-Paley
decomposition. Several related embedding relations and inequalities
will also be given here. For more details the reader is referred to
the
books~\cite{JbLo76},~\cite{Ca04},~\cite{Ch98},~\cite{Pe76},~\cite{Tr83}
for a comprehensive treatment.
\subsection{Littlewood-Paley Decomposition:} Let us start with the Littlewood-Paley decomposition
in $\mathbb{R}^{3}$. To this end, we take an arbitrary function
$\psi$ in the Schwartz class $\mathcal{S}(\mathbb{R}^{3}) $ whose
Fourier transform $\hat{\psi}$ is such that
\begin{eqnarray}
supp\ \hat{\psi} \subset \{\xi, \frac{1}{2}\leq \mid\xi\mid\leq 2\},
\end{eqnarray}
and
\begin{align} \forall\xi\neq 0, \quad\sum_{j \in\mathbb{Z}}
\hat{\psi}(\frac{\xi}{2^{j}}) = 1.\nonumber
\end{align} Let us define
$\varphi$ by
\begin{align}
\hat{\varphi}(\xi) = 1 - \sum_{j\geq
0}\hat{\psi}(\frac{\xi}{2^{j}}),\nonumber
\end{align}
and hence
\begin{eqnarray}
supp\ \hat{\varphi} \subset \{\xi, \mid\xi\mid\leq 1\}.
\end{eqnarray}
For $j\in \mathbb{Z}$, we write $\varphi_{j}(x) =
2^{3j}\varphi(2^{j}x)$. We denote by $S_{j}$ and $\triangle_{j}$,
the convolution operators with $\varphi_{j}$ and $\psi_{j}$
respectively. Hence
\begin{align} S_{j}(f) =
f\star\varphi_{j},\nonumber
\end{align}
and
\begin{align}
\triangle_{j}f = \psi_{j}\star f,\quad\text{where}\ \psi_{j}(x) =
2^{3j}\psi(2^{j}x).\nonumber
\end{align}
 Then
\begin{align}
 S_{j} = \sum_{p<j} \triangle_{p}\quad\text{and}\quad
 I = \sum_{j\in \mathbb{Z}}\ \triangle_{j}.\nonumber
\end{align}
The dyadic decomposition
\begin{eqnarray}\label{e3.3}
 u = \sum_{j\in \mathbb{Z}}\ \triangle_{j} u,
\end{eqnarray}
is called the homogeneous Littlewood-Paley decomposition of $u$ and
converges only in the quotient space
$\mathcal{S}^{'}/_{\mathcal{P}}$ where $\mathcal{S}^{'}$ is the
space of tempered distributions and $\mathcal{P}$ is the space of
polynomials. Now let us mention here the following
quasi-orthogonality properties of the dyadic
decomposition~\cite{Ch98} (proposition 2.1.1):
\begin{eqnarray}\label{e3.4}
\triangle_{p}\triangle_{q}u = 0 \quad\text{if}\quad\mid p-q \mid
\geq 2,
\end{eqnarray}
\begin{eqnarray}\label{e3.5}
\triangle_{p}(S_{q-2}u\triangle_{q}u) = 0 \quad\text{if}\quad\mid
p-q \mid \geq 4.
\end{eqnarray}
\subsection{Besov Spaces:}
Let $0<p, q\leq\infty$ and $s\in \mathbb{R}$. Then a tempered
distribution $f$ belongs to the homogeneous Besov space
$\dot{B}_{p}^{s,q}$ if and only if
\begin{eqnarray}\label{e3.6}
\Big(\sum_{j\in\mathbb{Z}}\
2^{jsq}\parallel\triangle_{j}f\parallel_{p}^{q}\Big)^{\frac{1}{q}} <
\infty
\end{eqnarray}
and $f = \sum_{j\in Z}\triangle_{j}f$ in $
\mathcal{S}^{'}/_{\mathcal{P}_{m}}$ where $\mathcal{P}_{m}$ is the
space of polynomials of degree $\leq m$ and $m = [s-\frac{d}{p}]$,
the integer part of $s-\frac{d}{p}$.

Besov space is a quasi-Banach space~\cite{RuSi96}. Here we recall
the following standard embedding
rules~\cite{Tr83} (chapter 2.7):\\
If $s_{1} > s_{2}$ and $p_{2}\geq p_{1}\geq 1$ such that
$s_{1}-\frac{d}{p_{1}} = s_{2}-\frac{d}{p_{2}}$, then
\begin{eqnarray}
\dot{B}_{p_{1}}^{s_{1},q}\hookrightarrow \dot{B}_{p_{2}}^{s_{2},q}.
\end{eqnarray}
Moreover if $q_{1} < q_{2}$ then
\begin{eqnarray}
\dot{B}_{p}^{s,q_{1}}\hookrightarrow \dot{B}_{p}^{s,q_{2}}.
\end{eqnarray}
The above mentioned embeddings are also valid for inhomogeneous
Besov spaces. For more embedding theorems and their proofs we refer
the readers to~\cite{Pe76} and~\cite{Tr83}.

 Next let us recall the following result from Chapter 3 in
 Triebel~\cite{Tr83}:

\begin{lem}
Let $1\leq p, q\leq\infty$ and $s<0$. Then $\forall f\in
\dot{B}_{p}^{s,q}$ we have,
\begin{eqnarray}\label{e3.9}
\Big(\sum_{j\in\mathbb{Z}}\
2^{jsq}\parallel\triangle_{j}f\parallel_{p}^{q}\Big)^{\frac{1}{q}} <
\infty \ \Leftrightarrow \ \Big(\sum_{j\in\mathbb{Z}}\
2^{jsq}\parallel S_{j}f\parallel_{p}^{q}\Big)^{\frac{1}{q}} <
\infty.
\end{eqnarray}
\end{lem}

Now we recall the following versions of Bernstein inequalities
(chapter 3 in~\cite{Lr02}):

\begin{lem}
Let $1\leq p\leq\infty$. Then there exist constants $C_{0},
C_{1}, C_{2} > 0$ such that \\
\textbf{(a)} If $f$ has its frequency in a ball $\mathbb{B}(0,
\lambda)$ $(supp\ \mathcal{F}(f) \subset \mathbb{B}(0, \lambda))$
then
\begin{eqnarray}\label{e3.10}
\parallel (-\triangle)^{\frac{s}{2}} f\parallel_{p} \ \leq C_{0} \ \lambda^{\mid
s\mid}\parallel f\parallel_{p}.
\end{eqnarray}
\textbf{(b)} If $f$ has its frequency in an annulus $\mathbb{C}(0,
A\lambda, B\lambda)$ \\ $(supp\ \mathcal{F}(f) \subset \{\xi,
A\lambda \leq \mid\xi\mid \leq B\lambda\} )$ then
\begin{eqnarray}\label{e3.11}
C_{1} \ \lambda^{\mid s\mid}\parallel f\parallel_{p} \ \leq \
\parallel (-\triangle)^{\frac{s}{2}} f\parallel_{p} \ \leq C_{2} \ \lambda^{\mid
s\mid}\parallel f\parallel_{p}.
\end{eqnarray}
\end{lem}

Now let us state here the modified Poincar\'{e} type inequality
given by Planchon~\cite{Pl00}.

\begin{lem}
Let $f\in \mathcal{S}$, the Schwartz space, whose fourier transform
is supported outside the ball $\mathbb{B}(0,1)$. Then for $p\geq 2$,
\begin{eqnarray}\label{e3.12}
\int\mid f\mid^{p}\d x\ \leq \ C_{p}\int\mid\nabla f\mid^{2} \mid
f\mid ^{p-2}\d x.
\end{eqnarray}
\end{lem}

\subsection{The Paraproduct rule:}
Another important tool in Littlewood-Paley analysis is the
paraproduct operator introduced by J. M. Bony~\cite{Bo81}. The idea
of the paraproduct enables us to define a new product between
distributions which turns out to be continuous in many functional
spaces where the pointwise multiplication does not make sense. This
is a
powerful tool for the analysis of nonlinear partial differential equations. \\
Let $f, g\in \mathcal{S}^{'}$. Then using the formal
Littlewood-Paley decomposition,
\begin{eqnarray}
 f = \sum_{j\in \mathbb{Z}} \triangle_{j} f, \qquad g = \sum_{j\in \mathbb{Z}} \triangle_{j}
 g.\nonumber
\end{eqnarray}
Hence
\begin{eqnarray}
fg &=& \sum_{j, l}\triangle_{j} f  \triangle_{l}
 g\nonumber\\
 &=& \sum_j\sum_{l<j-2}\triangle_{j} f  \triangle_{l}
 g + \sum_j\sum_{l>j+2}\triangle_{j} f  \triangle_{l}
 g + \sum_j\sum_{\mid l-j\mid\leq 2}\triangle_{j} f  \triangle_{l}
 g \nonumber\\
 &=& \sum_j\sum_{l<j-2}\triangle_{j} f  \triangle_{l}
 g + \sum_l\sum_{j<l-2}\triangle_{j} f  \triangle_{l}
 g + \sum_j\sum_{\mid l-j\mid\leq 2}\triangle_{j} f  \triangle_{l}
 g \nonumber\\
&=& \sum_{j}\triangle_{j}f S_{j-2}g + \sum_{j}\triangle_{j}g
S_{j-2}f + \sum_{\mid l-j\mid \leq 2}\triangle_{j}f
\triangle_{l}g\nonumber.
\end{eqnarray}
In other words, the product of two tempered distributions is
decomposed into two homogeneous paraproducts, respectively
\begin{eqnarray}
\dot\pi(f,g) = \sum_{j}\triangle_{j}f S_{j-2}g \qquad \text{and}
\qquad \dot\pi(g,f) = \sum_{j}\triangle_{j}g S_{j-2}f \nonumber,
\end{eqnarray}
plus a remainder
\begin{eqnarray}
R(f,g) = \sum_{\mid l-j\mid \leq 2}\triangle_{j}f
\triangle_{l}g\nonumber.
\end{eqnarray}

$\dot\pi$ is called the homogeneous paraproduct operator and the
convergence of the above series holds true in the quotient space
$\mathcal{S}^{'}/_{\mathcal{P}}$. Finally, using the
quasi-orthogonality properties from \eqref{e3.4} and \eqref{e3.5}
and after neglecting some non-diagonal terms for simplicity (since
the contributions from these non-diagonal terms are taken care of by
the terms which are being considered and hence negligible and also
this does not affect the convergence of the paraproducts~\cite{Ca95,
CaMe95}), we obtain
\begin{eqnarray}\label{e3.13}
\triangle_{j}(fg) = \triangle_{j}f S_{j-2}g + \triangle_{j}g
S_{j-2}f +  \triangle_{j}\big(\sum_{k\geq j}\triangle_{k}f
\triangle_{k}g\big).
\end{eqnarray}
We refer the
readers~\cite{Ca95},~\cite{Ch98},~\cite{Me81},~\cite{RuSi96} for
extensive studies on paraproducts.


\section{Main Results}

\begin{thm}[Local Lyapunov Property in $\mathrm{L}^p$]
Let $m\geq 3$, $2\leq p <\infty$. Let $\omega$ be the solution of
the vorticity equation \eqref{e1.1} such that
\begin{align}
& u\in \mathrm{C}{([0,T];\mathrm{L}^{m}\cap \mathrm{L}^{p})},
\qquad\nabla u\in \mathrm{L}^{1}_{loc}((0, T); \mathrm{L}^{p}),
\nonumber
\end{align}
and
\begin{align}
&\omega\in \mathrm{C}{([0,T];\mathrm{L}^{m}\cap
\mathrm{L}^{p})},\qquad\nabla\omega\in \mathrm{L}^{1}_{loc}((0, T);
\mathrm{L}^{p}),\quad\text{for} \ 0 < T \leq\infty.\nonumber
\end{align}
Then
\begin{eqnarray}\label{e4.1}
\partial_{t}\parallel\omega(t)\parallel_{p}^{p} &\leq&  -C(\nu - K \parallel
u(t)\parallel_{m})Q_{p}(\omega(t)) , \qquad 0 < t < T,
\end{eqnarray}
where $K$ denotes a positive constant depending upon $m$ and $p$.
\end{thm}
This implies for small $\mathrm{L}^m$-norm
$t\rightarrow\parallel\omega(t)\parallel_{\mathrm{L}^p}$ is a
Lyapunov function.

\begin{proof}
Consider,
\begin{eqnarray}
\partial_{t}\parallel\omega\parallel_{p}^{p} \nonumber &=&
\frac{\partial}{\partial t}\int\mid\omega\mid^{p}\d x \ = \
\frac{\partial}{\partial t}\int\mid\omega^{2}\mid^{p/2}\d x\nonumber
\\ &=& p \int \mid\omega\mid^{p-2}\omega\cdot\frac{\partial \omega}{\partial
t}\d x \ = \  p\langle\mid\omega\mid^{p-2}\omega,
\partial_{t}\omega\rangle\nonumber
\\ &=& p\langle\mid\omega\mid^{p-2}\omega,\ \nu\triangle\omega -
u\cdot\nabla\omega + \omega\cdot\nabla u\rangle\nonumber
\\ &=& \nu p \langle\mid\omega\mid^{p-2}\omega, \
\triangle\omega\rangle \ - \ p \langle\mid\omega\mid^{p-2}\omega, \
u\cdot\nabla\omega\rangle\nonumber\\
&&\quad + \ p \langle\mid\omega\mid^{p-2}\omega, \ \omega\cdot\nabla
u\rangle.\label{e4.2}
\end{eqnarray}
Using Lemma 2.5 on the first term of the right hand side, we have
from \eqref{e4.2}
\begin{eqnarray}\label{e4.3}
\partial_{t}\parallel\omega\parallel_{p}^{p}\ &\leq&
-C \nu Q_{p}(\omega) - \ p \langle\mid\omega\mid^{p-2}\omega, \
u\cdot\nabla\omega\rangle \ + \ p \langle\mid\omega\mid^{p-2}\omega,
\ \omega\cdot\nabla u\rangle\nonumber\\
\end{eqnarray}
Now we need to estimate the second and the third terms of the right
hand side of the equation \eqref{e4.3}.

Using the fact that $\Div u = 0$, we have
\begin{eqnarray}\label{e4.4}
u\cdot\nabla\omega \ = \ u_{i}\frac{\partial\omega_{j}}{\partial
x_{i}} \ = \ \frac{\partial}{\partial x_{i}}(u_{i}\omega_{j}) \ - \
\omega_{j}\frac{\partial u_{i}}{\partial x_{i}} \ =
\nabla\cdot(u\otimes\omega),
\end{eqnarray}
where $\otimes$ represents the tensor product.\\
Then
\begin{eqnarray}
\mid\langle\mid\omega\mid^{p-2}\omega, \
u\cdot\nabla\omega\rangle\mid &=&
\mid\langle\mid\omega\mid^{p-2}\omega, \
\nabla\cdot(u\otimes\omega)\rangle\mid\nonumber
\\ &=& \mid\langle\nabla(\mid\omega\mid^{p-2}\omega), \
 u\otimes\omega\rangle\mid\nonumber
\\ &\leq & \langle\mid\nabla(\mid\omega\mid^{p-2}\omega)\mid, \
 \mid u \otimes\omega\mid\rangle.\label{e4.5}
\end{eqnarray}
Notice that $ \mid\nabla\mid\omega\mid^{p-2}\omega\mid
 \ \leq\ C\mid\omega\mid^{p-2}\mid\nabla\omega\mid.$
 Hence using this and H\"{o}lder's inequality in \eqref{e4.5} we have
\begin{eqnarray}
\mid\langle\mid\omega\mid^{p-2}\omega, \
u\cdot\nabla\omega\rangle\mid &\leq&
\langle\mid\omega\mid^{p-2}\mid\nabla\omega\mid, \ \mid u
 \otimes\omega\mid\rangle\nonumber\\
 &\leq &
 \parallel\ \mid\omega\mid^{p-2}\mid\nabla\omega\mid \ \parallel_{q} \
 \parallel u\otimes\omega\parallel_{q'},\quad \text{where}\
 \frac{1}{q}+\frac{1}{q'}=1.\nonumber\\ \label{e4.6}
\end{eqnarray}
Now
\begin{align} \parallel \
\mid\omega\mid^{p-2}\mid\nabla\omega\mid \
\parallel_{q}^{q}
&=\int \mid\omega\mid^{q(p-2)} \ \mid\nabla\omega\mid^{q}\d
x\nonumber
\\&=\int \mid\omega\mid^{q(p-2)/2} \ \big(\ \mid\omega\mid^{p-2} \
\mid\nabla\omega\mid^{2}\ \big)^{q/2}\d x.\nonumber
\end{align}
Since $\frac{2-q}{2}+\frac{q}{2}=1$, H\"{o}lder inequality yields
\begin{align}
&\parallel \ \mid\omega\mid^{p-2}\mid\nabla\omega\mid \
\parallel_{q}^{q}\nonumber
\\&\leq\Big[\int \big(\ \mid\omega\mid^{q(p-2)/2}\
\big)^{2/(2-q)}\d x \Big]^\frac{2-q}{2} \ \Big[\int \Big\{\big(\
\mid\omega\mid^{p-2} \ \mid\nabla\omega\mid^{2}\big)^{q/2}
\Big\}^{2/q}\d x \Big ]^\frac{q}{2} \nonumber
\\&=\Big[\int\mid\omega\mid^{q(p-2)/(2-q)}\d x\Big]^{(2-q)/2} \ \Big[
\int\mid\omega\mid^{p-2}\ \mid\nabla\omega\mid^{2}\d
x\Big]^{q/2}\nonumber
\\&=\ \parallel\omega\parallel_{r}^{q(p-2)/2} \ Q_{p}(\omega)^{q/2},
\quad \text{where} \ r = \frac{q(p-2)}{(2-q)}. \nonumber
\end{align}
Hence
\begin{align}\label{e4.7}
\parallel \ \mid\omega\mid^{p-2}\mid\nabla\omega\mid \
\parallel_{q}\ \leq\ \parallel\omega\parallel_{r}^{(p-2)/2} \
Q_{p}(\omega)^{1/2}.
\end{align}
Again by H\"{o}lder,
\begin{align}\label{e4.8}
\parallel u\otimes\omega\parallel_{q'}\ \leq\ C\parallel
u\parallel_{m} \ \parallel\omega\parallel_{r},\quad \text{since} \
\frac{1}{q'}=\frac{1}{m}+\frac{1}{r}.
\end{align}
Now from the relations
\begin{eqnarray}
\frac{1}{q} + \frac{1}{q'} = 1 ,\ r = \frac{q(p-2)}{(2-q)}\
\text{and} \ \frac{1}{q'} = \frac{1}{m} + \frac{1}{r}, \nonumber
\end{eqnarray}
we find that
\begin{eqnarray}\label{e4.9}
\qquad\ r = \frac{mp}{(m-2)}.
\end{eqnarray}
Using equations \eqref{e4.7} and \eqref{e4.8} in \eqref{e4.6} we
have
\begin{eqnarray}
\mid\langle\mid\omega\mid^{p-2}\omega, \
u\cdot\nabla\omega\rangle\mid &\leq &
C\parallel\omega\parallel_{r}^{(p-2)/2} \ Q_{p}(\omega)^{1/2}
\parallel u\parallel_{m}\ \parallel\omega\parallel_{r}\nonumber
\\&=& C\parallel u\parallel_{m}\ \parallel\omega\parallel_{r}^{p/2}\
Q_{p}(\omega)^{1/2}.\nonumber
\end{eqnarray}
Applying the Lemma 2.6 in the above equation we obtain
\begin{eqnarray}\label{e4.10}
\mid\langle\mid\omega\mid^{p-2}\omega, \
u\cdot\nabla\omega\rangle\mid &\leq& C\parallel u\parallel_{m}\
Q_{p}(\omega).
\end{eqnarray}
The third term in the equation \eqref{e4.3} can be estimated by
using the fact that $\Div \omega = 0$ along with the similar kind of
techniques taken to estimate the second term.

Thus we get
\begin{eqnarray}\label{e4.11}
\mid\langle\mid\omega\mid^{p-2}\omega,\ \omega\cdot\nabla\ u
\rangle\mid &\leq & C\parallel u\parallel_{m}\ Q_{p}(\omega).
\end{eqnarray}
Combining \eqref{e4.10} and \eqref{e4.11} with \eqref{e4.3} we get
the desired result \eqref{e4.1}.
\end{proof}

\begin{thm}[Local Lyapunov Property in Besov Spaces]
 Let the initial data $\omega_{0}$ for the $3$-D vorticity equation be
in $\dot{B}_{p}^{s,q}$ where $ s= \frac{3}{p}-1$, $p, q \geq 2$, and
$\frac{3}{p}+\frac{2}{q}>1$. Then there exist small constants
$\varepsilon_{1}> 0$ and $\varepsilon_{2}> 0$ such that if the
velocity field satisfies $sup_{t}\parallel u(t)
\parallel_{\dot{B}_{\infty}^{-1, \infty}} < \varepsilon_{1}$ and the
vorticity field satisfies $sup_{t}\parallel \omega
(t)\parallel_{\dot{B}_{\infty}^{-2, \infty}} < \varepsilon_{2}$,
then $ t\rightarrow\parallel\omega(t)\parallel_{\dot{B}_{p}^{s,q}}$
is a Lyapunov function.
\end{thm}

\begin{proof}
Let us consider
\begin{align}
F(u,w) = u\cdot\nabla\omega - \omega\cdot\nabla u.\nonumber
\end{align}
Multiply the equation \eqref{e1.1} by $\triangle_{j}$ to get,
\begin{align}\label{e4.12}
\partial_{t}(\triangle_{j}\omega)\ - \nu\triangle(\triangle_{j}\omega)\
+ \triangle_{j}(F(u,w)) = 0.
\end{align}
Now,
\begin{align}
\partial_{t}\parallel\triangle_{j}\omega\parallel_{p}^{p}&=
\frac{\partial}{\partial t}\int\mid\triangle_{j}\omega\mid^{p}\d x \
= \ \frac{\partial}{\partial
t}\int\mid(\triangle_{j}\omega)^{2}\mid^{p/2}\d x\nonumber
\\ &=p \int \mid\triangle_{j}\omega\mid^{p-2}\triangle_{j}\omega\cdot\frac{\partial(\triangle_{j}\omega)}{\partial
t}\d x\nonumber\\
&=p\langle\mid\triangle_{j}\omega\mid^{p-2}\triangle_{j}\omega,
\partial_{t}(\triangle_{j}\omega)\rangle.\nonumber
\end{align}
Hence using \eqref{e4.12} we have from the above equation
\begin{align}
\partial_{t}\parallel\triangle_{j}\omega\parallel_{p}^{p}
&=p\langle\mid\triangle_{j}\omega\mid^{p-2}\triangle_{j}\omega,\
\nu\triangle(\triangle_{j}\omega) - \triangle_{j}(F(u,w))
\rangle\nonumber\\
&=\nu p \langle\mid\triangle_{j}\omega\mid^{p-2}\triangle_{j}\omega,
\ \triangle(\triangle_{j}\omega)\rangle \nonumber\\
&\qquad -p
\langle\mid\triangle_{j}\omega\mid^{p-2}\triangle_{j}\omega,
\triangle_{j}(F(u,w)) \rangle.\nonumber
\end{align}
Applying the Lemma 2.5 on the first term on the right hand side of
the above equation we obtain
\begin{align}
\partial_{t}\parallel\triangle_{j}\omega\parallel_{p}^{p}
&\leq -\nu
p\int\mid\triangle_{j}\omega\mid^{p-2}\mid\nabla\triangle_{j}\omega\mid^{2}dx
\nonumber\\
&\qquad -p
\langle\mid\triangle_{j}\omega\mid^{p-2}\triangle_{j}\omega,
\triangle_{j}(F(u,w)) \rangle.\nonumber
\end{align}
Hence,
\begin{align}
\partial_{t}\parallel\triangle_{j}\omega\parallel_{p}^{p} +
&\nu
p\int\mid\triangle_{j}\omega\mid^{p-2}\mid\nabla\triangle_{j}\omega\mid^{2}\d
x\nonumber\\
&\leq -p \int\mid\triangle_{j}\omega\mid^{p-2}\triangle_{j}\omega
\triangle_{j}(F(u,w))\d x, \nonumber
\end{align}
which is equivalent of considering the equation
\begin{align}
\frac{d}{dt}\parallel\triangle_{j}\omega\parallel_{p}^{p} + & \nu
p\int\mid\triangle_{j}\omega\mid^{p-2}\mid\nabla\triangle_{j}\omega\mid^{2}\d
x\nonumber\\
&\leq p
\int\mid\triangle_{j}\omega\mid^{p-1}\mid\triangle_{j}(F(u,w))\mid
dx. \nonumber
\end{align}
Using Lemma 3.3 we replace the second term to get
\begin{align}\label{e4.13}
&\frac{d}{dt}\parallel\triangle_{j}\omega\parallel_{p}^{p} + \
\tilde{C}_p\nu p\ 2^{2j}
\parallel\triangle_{j}\omega\parallel_{p}^{p} \ \leq p
\int\mid\triangle_{j}\omega\mid^{p-1}\mid\triangle_{j}(F(u,w))\mid
\d x,\\
&\text{where $\tilde{C}_p$ is positive constant depending on
$p$}.\nonumber
\end{align}
Now,
\begin{align}
\mid\triangle_{j}(F(u,w))\mid\nonumber &=
\mid\triangle_{j}(u\cdot\nabla\omega - \omega\cdot\nabla u)\mid
\nonumber
\\ &\leq\mid\triangle_{j}(u\cdot\nabla\omega)\mid + \mid\triangle_{j}(\omega\cdot\nabla
u)\mid.\nonumber
\end{align}
Moreover
\begin{align}
u\cdot\nabla\omega \ = \ u_{i}\frac{\partial\omega_{j}}{\partial
x_{i}} \ = \ \frac{\partial}{\partial x_{i}}(u_{i}\omega_{j}) \ - \
\omega_{j}\frac{\partial u_{i}}{\partial x_{i}} \ =
\nabla\cdot(u\otimes\omega), \quad \text{since}\ \Div u =
0\nonumber,
\end{align}
and similarly $\omega\cdot\nabla u = \nabla\cdot(\omega\otimes u),$
where $\otimes$ represents the usual tensor product. \\
Since the terms $\nabla\cdot(u\otimes\omega)$ and
$\nabla\cdot(\omega\otimes u)$ behave in similar fashion, we have
from equation \eqref{e4.13}
\begin{align}\label{e4.14}
\frac{d}{dt}\parallel\triangle_{j}\omega\parallel_{p}^{p} + \
\tilde{C}_p\nu p\ 2^{2j}
\parallel\triangle_{j}\omega\parallel_{p}^{p} \ \leq\ 2p
\int\mid\triangle_{j}\omega\mid^{p-1}\mid\triangle_{j}\nabla\cdot(u\otimes\omega)\mid
\d x.
\end{align}
Now using the paraproduct rule \eqref{e3.13}, we have
\begin{align}
\triangle_{j}\nabla\cdot(u\otimes\omega)\nonumber &=
\nabla\triangle_{j}(u\otimes\omega)\nonumber
\\ &= \nabla\big(\triangle_{j}u \ S_{j-2}\omega\big) +
\nabla\big(\triangle_{j}\omega \ S_{j-2}u\big) + \nabla\big(
\triangle_{j}\big(\sum_{k\geq j}\triangle_{k}u\
\triangle_{k}\omega\big)\big).\label{e4.15}
\end{align}
Using \eqref{e4.15} in \eqref{e4.14} we obtain,
\begin{align}
&\frac{d}{dt}\parallel\triangle_{j}\omega\parallel_{p}^{p} + \
\tilde{C}_p
\nu p\ 2^{2j} \parallel\triangle_{j}\omega\parallel_{p}^{p}\nonumber\\
&\quad\leq 2p
\int\mid\triangle_{j}\omega\mid^{p-1}\mid\nabla\big(\triangle_{j}u \
S_{j-2}\omega\big)\mid \d x \nonumber
\\ &\quad\quad + 2p \int\mid\triangle_{j}\omega\mid^{p-1}\mid\nabla\big(\triangle_{j}\omega
\ S_{j-2}u\big)\mid \d x \nonumber
\\ &\quad\quad + 2p \int\mid\triangle_{j}\omega\mid^{p-1}\mid\nabla\big(\triangle_{j}\big(\sum_{k\geq
j}\triangle_{k}u\ \triangle_{k}\omega\big)\big)\mid \d
x.\label{e4.16}
\end{align}
We need to estimate each of the terms on the right hand side of
\eqref{e4.16} separately.

First consider the term
\begin{align}
\int\mid\triangle_{j}\omega\mid^{p-1}\mid\nabla\big(\triangle_{j}\omega
\ S_{j-2}u\big)\mid \d x,\nonumber
\end{align}
 and apply H\"{o}lder's Inequality to get,
\begin{align}
\int\mid\triangle_{j}\omega\mid^{p-1}\mid\nabla\big(\triangle_{j}\omega
\ S_{j-2}u\big)\mid \d x\leq \
\parallel\triangle_{j}\omega\parallel_{p}^{p-1} \
\parallel\nabla\big(\triangle_{j}\omega
\ S_{j-2}u\big)\parallel_{p}.\nonumber
\end{align}
With the help of Lemma 3.2 we obtain
\begin{align}
&\int\mid\triangle_{j}\omega\mid^{p-1}\mid\nabla\big(\triangle_{j}\omega
\ S_{j-2}u\big)\mid \d x\nonumber\\
&\quad\leq \ C_1 \parallel\triangle_{j}\omega\parallel_{p}^{p-1} \
2^{j}
\parallel\triangle_{j}\omega\ S_{j-2}u\parallel_{p}\nonumber\\
&\quad=\ C_1 \parallel\triangle_{j}\omega\parallel_{p}^{p-1} \ 2^{j}
\parallel\big(2^{j}\triangle_{j}\omega\big)\big(2^{-j}S_{j-2}u\big)\parallel_{p}\nonumber
\\ &\quad\leq\ C_1 \parallel\triangle_{j}\omega\parallel_{p}^{p-1} \ 2^{j}
\parallel 2^{j}\triangle_{j}\omega\parallel_{p} \ \sup_{j}\big(2^{-j}\parallel
S_{j-2}u\parallel_{\infty}\big)\nonumber
\\ &\quad=\ C_1\ 2^{2j}\parallel\triangle_{j}\omega\parallel_{p}^{p}\ \sup_{j}\big(2^{-j}\parallel
S_{j-2}u\parallel_{\infty}\big),\label{e4.17}\\
&\text{where $C_1$ is positive constant.}\nonumber
\end{align}

 Now from Lemma 3.1, for $s = -1$
and $p = q = \infty$, we have
\begin{align}
&\qquad 2^{-j}\parallel\triangle_{j}
u\parallel_{\mathrm{L}^{\infty}}\ \in l^{\infty}\ \Leftrightarrow \
2^{-j}\parallel S_{j}
u\parallel_{\mathrm{L}^{\infty}}\ \in l^{\infty}\nonumber\\
&\Rightarrow \ \sup_{j} 2^{-j}\parallel\triangle_{j}
u\parallel_{\infty}\ \Leftrightarrow \ \sup_{j} 2^{-j}\parallel
S_{j} u\parallel_{\infty}\nonumber\\
&\Rightarrow \ \parallel u(x,t)\parallel_{\dot{B}_{\infty}^{-1,
\infty}}\ \Leftrightarrow \ \sup_{j} 2^{-j}\parallel S_{j}
u\parallel_{\infty}.
\end{align}
Then using the conditions assumed in the theorem, we get,
\begin{align}
&\qquad\parallel u(x,t)\parallel_{\dot{B}_{\infty}^{-1, \infty}} \
\leq \ \sup_{t}\parallel u(x,t)\parallel_{\dot{B}_{\infty}^{-1,
\infty}}\ \leq\varepsilon_{1},\nonumber\\
&\Rightarrow \ \sup_{j} 2^{-j}\parallel S_{j} u\parallel_{\infty}\
\leq\varepsilon_{1}.
\end{align}
So finally \eqref{e4.17} yields
\begin{align}\label{e4.20}
\int\mid\triangle_{j}\omega\mid^{p-1}\mid\nabla\big(\triangle_{j}\omega
\ S_{j-2}u\big)\mid dx \ \leq C_1\varepsilon_{1}
2^{2j}\parallel\triangle_{j}\omega\parallel_{p}^{p}.
\end{align}
Now let us consider the term:
\begin{align}
\int\mid\triangle_{j}\omega\mid^{p-1}\mid\nabla\big(\triangle_{j}u \
S_{j-2}\omega\big)\mid \d x.\nonumber
\end{align}
As before H\"{o}lder's Inequality and Lemma 3.2 yield
\begin{align}
&\int\mid\triangle_{j}\omega\mid^{p-1}\mid\nabla\big(\triangle_{j}u
\ S_{j-2}\omega\big)\mid \d x\nonumber\\
&\quad\leq \ \parallel\triangle_{j}\omega\parallel_{p}^{p-1} \
\parallel\nabla\big(\triangle_{j} u
\ S_{j-2}\omega\big)\parallel_{p}\nonumber\\
&\quad\leq \ C_2 \parallel\triangle_{j}\omega\parallel_{p}^{p-1} \
2^{j}
\parallel\triangle_{j} u\ S_{j-2}\omega\parallel_{p}\nonumber\\
&\quad=\ C_2 \parallel\triangle_{j}\omega\parallel_{p}^{p-1} \ 2^{j}
\parallel\big(2^{2j}\triangle_{j} u\big)\big(2^{-2j}S_{j-2}\omega\big)\parallel_{p}\nonumber\\
&\quad\leq\ C_2 \parallel\triangle_{j}\omega\parallel_{p}^{p-1} \
2^{j}
\parallel 2^{2j}\triangle_{j} u\parallel_{p} \ \sup_{j}\big(2^{-2j}\parallel
S_{j-2}\omega\parallel_{\infty}\big),\label{e4.21}\\
&\text{where $C_2$ is positive constant.}\nonumber
\end{align}
From Lemma 3.2, equation \eqref{e3.11}, we obtain
\begin{align}\label{e4.22}
2^{j}\parallel\triangle_{j} u\parallel_{p}\ \leq \
\parallel\nabla\triangle_{j} u\parallel_{p}.
\end{align}
The above equation and \eqref{e2.5} in Lemma 2.7 yield
\begin{align}\label{e4.24}
\parallel\triangle_{j} u\parallel_{p}\ \leq \ 2^{-j}\parallel
\triangle_{j}\omega\parallel_{p}.
\end{align}
Now applying Lemma 3.1, for $s = -2$ and $p = q = \infty$ and
proceeding as before we obtain
\begin{align}\label{e4.25}
\sup_{j} 2^{-j}\parallel S_{j}\omega\parallel_{\infty}\
\leq\varepsilon_{2}.
\end{align}
Using \eqref{e4.24} and \eqref{e4.25} in \eqref{e4.21} we have
\begin{align}\label{e4.26}
\int\mid\triangle_{j}\omega\mid^{p-1}\mid\nabla\big(\triangle_{j}u \
S_{j-2}\omega\big)\mid dx \leq C_2\varepsilon_{2}
2^{2j}\parallel\triangle_{j}\omega\parallel_{p}^{p}.
\end{align}
Next we estimate the last term
\begin{align}
&\int\mid\triangle_{j}\omega\mid^{p-1}\mid\nabla\big(\triangle_{j}\big(\sum_{k\geq
j}\triangle_{k}u\ \triangle_{k}\omega\big)\big)\mid \d x\nonumber\\
&\quad\leq \ \parallel\triangle_{j}\omega\parallel_{p}^{p-1}
\parallel\nabla\big(\triangle_{j}\big(\sum_{k\geq j}\triangle_{k}u\
\triangle_{k}\omega\big)\big)\parallel_{p}\nonumber\\
&\quad\leq\ C_3 2^{j}\parallel\triangle_{j}\omega\parallel_{p}^{p-1}
\parallel\triangle_{j}\big(\sum_{k\geq j}\triangle_{k}u\
\triangle_{k}\omega\big)\parallel_{p}\nonumber\\
&\text{where $C_3$ is positive constant.}\nonumber
\end{align}
Using Young's Inequality as in~\cite{CaMe95}, we have
\begin{align}
&\int\mid\triangle_{j}\omega\mid^{p-1}\mid\nabla\big(\triangle_{j}\big(\sum_{k\geq
j}\triangle_{k}u\ \triangle_{k}\omega\big)\big)\mid \d
x\nonumber\\&\quad\leq\ C_p
2^{j}\parallel\triangle_{j}\omega\parallel_{p}^{p-1}
\big(\sum_{k\geq j}\parallel\triangle_{k}u\parallel_{p}\parallel
\triangle_{k}\omega\parallel_{p}\big),\label{e4.27}\\
&\text{where $C_p$ is positive constant depending on $p$.}\nonumber
\end{align}
Now
\begin{align}
\parallel\triangle_{j}\omega\parallel_{p}^{p-1} &=\ 2^{2j}
\big(2^{-2j}\parallel\triangle_{j}\omega\parallel_{p}\big)\
\parallel\triangle_{j}\omega\parallel_{p}^{p-2}\nonumber
\\ &\leq\ 2^{2j} \big(\sup_{j}
2^{-2j}\parallel\triangle_{j}\omega\parallel_{\infty}\big) \
\parallel\triangle_{j}\omega\parallel_{p}^{p-2}\nonumber
\\ &=\ 2^{2j} \parallel\omega(x,t)\parallel_{\dot{B}_{\infty}^{-2,
\infty}}\
\parallel\triangle_{j}\omega\parallel_{p}^{p-2}\nonumber
\\ &\leq\ \varepsilon_{2}\
2^{2j}\parallel\triangle_{j}\omega\parallel_{p}^{p-2}.\label{e4.28}
\end{align}
Using \eqref{e4.24} and \eqref{e4.28} in \eqref{e4.27} we have,
\begin{align}\label{e4.29}
&\int\mid\triangle_{j}\omega\mid^{p-1}\mid\nabla\big(\triangle_{j}\big(\sum_{k\geq
j}\triangle_{k}u\ \triangle_{k}\omega\big)\big)\mid\d x\nonumber\\
&\quad\leq \ C_p \varepsilon_{2}\
2^{2j}\parallel\triangle_{j}\omega\parallel_{p}^{p-2}
\big(\sum_{k\geq j}\parallel\triangle_{k}
\omega\parallel_{p}^{2}\big).
\end{align}
Now combining all results from \eqref{e4.20}, \eqref{e4.26} and
\eqref{e4.29} and neglecting the constants $C_1, C_2, C_p,
\tilde{C}_p$, we obtain from \eqref{e4.16}
\begin{align}
\frac{d}{dt}\parallel\triangle_{j}\omega\parallel_{p}^{p} + \nu p\
2^{2j}
\parallel\triangle_{j}\omega\parallel_{p}^{p}\ &\leq\ 2p\varepsilon
_{1}\ 2^{2j}\parallel\triangle_{j}\omega\parallel_{p}^{p} +\
2p\varepsilon_{2}\
2^{2j}\parallel\triangle_{j}\omega\parallel_{p}^{p} \nonumber\\
&\quad +\ 2p\varepsilon_{2}\
2^{2j}\parallel\triangle_{j}\omega\parallel_{p}^{p-2}
\big(\sum_{k\geq j}\parallel\triangle_{k}
\omega\parallel_{p}^{2}\big).\nonumber
\end{align}
Simplifying we get,
\begin{align}\label{e4.30}
\frac{d}{dt}\parallel\triangle_{j}\omega\parallel_{p}^{2} +\ p
\big(\nu - 2\varepsilon_{1} -
2\varepsilon_{2}\big)2^{2j}\parallel\triangle_{j}\omega\parallel_{p}^{2}
\ &\leq\ 2p\varepsilon_{2}\ 2^{2j}\big(\sum_{k\geq
j}\parallel\triangle_{k} \omega\parallel_{p}^{2}\big).
\end{align}
The rest of the construction is motivated by~\cite{CaPl00}.
Multiplying both sides of \eqref{e4.30} by
$2^{jqs}\parallel\triangle_{j}\omega\parallel_{p}^{q-2}$, we get
\begin{align}
&\frac{d}{dt}\big(2^{jqs}\parallel\triangle_{j}\omega\parallel_{p}^{q}\big)
+\ p \big(\nu - 2\varepsilon_{1} -
2\varepsilon_{2}\big)2^{j(qs+2)}\parallel\triangle_{j}\omega\parallel_{p}^{q}\nonumber\\
&\quad\leq \ 2p\varepsilon_{2}\
2^{j(qs+2)}\parallel\triangle_{j}\omega\parallel_{p}^{q-2}\big(\sum_{k\geq
j}\parallel\triangle_{k} \omega\parallel_{p}^{2}\big).\nonumber
\end{align}
Let $ s= \frac{3}{p}-1$, $p, q \geq 2$, and
$\frac{3}{p}+\frac{2}{q}>1$. Then $r = \frac{2}{q} + s > 0$ or $qs +
2 = rq$. \\
Then
\begin{align}
&\frac{d}{dt}\big(2^{jqs}\parallel\triangle_{j}\omega\parallel_{p}^{q}\big)
+ \ p \big(\nu - 2\varepsilon_{1} -
2\varepsilon_{2}\big)2^{rqj}\parallel\triangle_{j}\omega\parallel_{p}^{q}\nonumber\\
&\quad\leq \ 2p\varepsilon_{2}\
2^{(q-2)rj}\parallel\triangle_{j}\omega\parallel_{p}^{q-2}
2^{2rj}\big(\sum_{k\geq j}\parallel\triangle_{k}
\omega\parallel_{p}^{2}\big).\label{e4.31}
\end{align}
Let $ f_{j} = 2^{js}\parallel\triangle_{j}\omega\parallel_{p}$ and
$g_{j} = 2^{jr}\parallel\triangle_{j}\omega\parallel_{p}$. Then
taking sum over $j$ of \eqref{e4.31} we have,
\begin{align}
\frac{d}{dt}\big(\sum_{j}f_{j}^{q}\big) + \ p \big(\nu -
2\varepsilon_{1} - 2\varepsilon_{2}\big)\sum_{j}g_{j}^{q}&\leq\
2p\varepsilon_{2} \sum_{j}g_{j}^{q-2} \ 2^{2rj}\big(\sum_{k\geq
j}\parallel\triangle_{k}
\omega\parallel_{p}^{2}\big)\nonumber\\
&=\ 2p\varepsilon_{2}\sum_{k=1}^{\infty}\sum_{j=1}^{k}\ g_{j}^{q-2}
\
2^{2rj}\parallel\triangle_{k} \omega\parallel_{p}^{2}\nonumber\\
&=\ 2p\varepsilon_{2}\sum_{k=1}^{\infty}\sum_{j=1}^{k}\ g_{j}^{q-2}
\ 2^{2rj} \ 2^{-2rk}\ g_{k}^{2}.\label{e4.32}
\end{align}
Let us consider $$\sum_{j=1}^{k}\ g_{j}^{q-2} \ 2^{2rj} = 2^{2rk}\
h_{k}^{q-2}.$$ Then it is clear that
\begin{align}\label{e4.33}
\sum_{k} h_{k}^{q} \ \leq \sum_{j} g_{j}^{q}.
\end{align}
So \eqref{e4.32} yields with the help of H\"{o}lder Inequality and
\eqref{e4.33}
\begin{align}
&\frac{d}{dt}\big(\sum_{j}f_{j}^{q}\big) + \ p \big(\nu -
2\varepsilon_{1} - 2\varepsilon_{2}\big)\sum_{j}g_{j}^{q}\nonumber\\
&\quad\leq\ 2p\varepsilon_{2}\sum_{k} h_{k}^{q-2}
g_{k}^{2}\nonumber\\
&\quad\leq\
2p\varepsilon_{2}\Big(\sum_{k}\big(h_{k}^{q-2}\big)^{\frac{q}{q-2}}\Big)^{\frac{q-2}{q}}
\
\Big(\sum_{k}\big(g_{k}^{2}\big)^{\frac{q}{2}}\Big)^{\frac{2}{q}},\quad\text{since}\
\ \frac{q-2}{q}+\frac{2}{q} =
 1,\nonumber\\
&\quad\leq\ 2p\varepsilon_{2} \big(\sum_{k}
g_{k}^{q}\big)^{\frac{q-2}{q}}\ \big(\sum_{k}
g_{k}^{q}\big)^\frac{2}{q}\nonumber\\
&\quad=\ 2p\varepsilon_{2}\sum_{k} g_{k}^{q}.
\end{align}
Hence,
\begin{align}
\frac{d}{dt}\big(\sum_{j}f_{j}^{q}\big) + \ p \big(\nu -
2\varepsilon_{1} - 4\varepsilon_{2}\big)\sum_{j}g_{j}^{q}\ \leq 0.
\end{align}
Using the definition of Besov Spaces in \eqref{e3.6}, we can write,
\begin{align}
\frac{d}{dt}\big(\parallel\omega(x,t)\parallel_{\dot{B}_{p}^{s,
q}}^{q}\big) + \ p \big(\nu - 2\varepsilon_{1} -
4\varepsilon_{2}\big)\parallel\omega(x,t)\parallel_{\dot{B}_{p}^{r,
q}}^{q} \ \leq 0.\nonumber
\end{align}
Hence $ t\rightarrow\parallel\omega(t)\parallel_{\dot{B}_{p}^{s,q}}$
is a Lyapunov function for small $\varepsilon_{1}$ and
$\varepsilon_{2}$ and comparatively large $\nu$.
\end{proof}

Now let us prove the dissipativity of the sum of the linear and
nonlinear operators of the vorticity equation in
$\mathrm{L}^{m}(\mathbb{R}^{m})$.

We write \eqref{e1.1} in the form $\partial_{t}\omega = \emph{A}
(u,\omega)$, where $\emph{A} : u, \omega \mapsto \emph{A} (u,\omega)
= \nu\triangle\omega - u\cdot\nabla\omega + \omega\cdot\nabla u$ is
a nonlinear operator. We know that $ G(\omega) =
\mid\omega\mid^{p-2}\omega$ \ is the duality map on $\mathrm{L}^{p}$
to $\mathrm{L}^{p^{\prime}}$. In Theorem 4.1 we proved that
\begin{align}\label{e4.36}
\langle\emph{A} (u,\omega), G(\omega)\rangle\ \leq\ -C(\nu - K
\parallel u(t)\parallel_{m})Q_{p}(\omega(t)).
\end{align}
Here we will prove a stronger property than \eqref{e4.36}.

\begin{thm}[Local Dissipativity in $\mathrm{L}^p$]
Let \ $m\geq 3$, $2\leq p <\infty$. Then if $(\omega -
\tilde{\omega})\in {\mathrm{L}^{1}}(\mathbb{R}^{m}) \cap
{\mathrm{L}^{r}}(\mathbb{R}^{m}), \ for \ r = \frac{mp}{m-2},$
\begin{align}
&\langle\emph{A}(u,\omega) - \emph{A}(v,\tilde{\omega}), G(\omega -
\tilde{\omega})\rangle\nonumber\\
&\quad\leq\ -C\big(\nu - K (\parallel u
\parallel _{m} +\parallel v \parallel _{m} + \parallel\omega
\parallel_{m}+\parallel \tilde{\omega} \parallel _{m})Q_{p}(\omega -
\tilde{\omega})\nonumber \\
&\quad\quad\ -K(\parallel\omega
\parallel_{m}+\parallel \tilde{\omega} \parallel
_{m})\parallel\omega -
\tilde{\omega}\parallel_{\mathrm{L}^{1}}Q_{p}(\omega -
\tilde{\omega})^{1/p'}\big),\label{e4.37}
\end{align}
where $\frac{1}{p} + \frac{1}{p'} = 1.$
\end{thm}

Hence, in the light of \eqref{e2.6}, we note that if $\omega$ and
$\tilde{\omega}$ are small in $\mathrm{L}^1\cap\mathrm{L}^m$, then
$$\langle\emph{A}(u,\omega) - \emph{A}(v,\tilde{\omega}), G(\omega -
\tilde{\omega})\rangle\ \leq\ 0,$$ which is a local dissipativity
property for $\emph{A}(\cdot, \cdot)$.

\begin{proof}
It is clear that
\begin{align}
&\langle\emph{A}(u,\omega) - \emph{A}(v,\tilde{\omega}),\ G(\omega -
\tilde{\omega})\rangle\nonumber\\
&\quad= \langle\nu\triangle\omega - u\cdot\nabla\omega +
\omega\cdot\nabla u -\nu\triangle\tilde{\omega} +
v\cdot\nabla\tilde{\omega} - \tilde{\omega}\cdot\nabla v, \ G(\omega
- \tilde{\omega})\rangle\nonumber\\
&\quad= \nu\langle\triangle(\omega-\tilde{\omega}),\ G(\omega -
\tilde{\omega})\rangle - \langle u\cdot\nabla\omega -
v\cdot\nabla\tilde{\omega},\ G(\omega -
\tilde{\omega})\rangle\nonumber\\
&\quad\quad + \langle\omega\cdot\nabla u - \tilde{\omega}\cdot\nabla
v,\ G(\omega - \tilde{\omega})\rangle.\label{e4.38}
\end{align}
According to Lemma 2.5
\begin{align}\label{e4.39}
\nu\langle\triangle(\omega-\tilde{\omega}),\ G(\omega -
\tilde{\omega})\rangle\ \leq\ -C\nu Q_{p}(\omega - \tilde{\omega}).
\end{align}
Now we need to estimate the second and third terms of the right hand
side of \eqref{e4.38}.\\
Notice that
\begin{align}
&\mid\langle u\cdot\nabla\omega - v\cdot\nabla\tilde{\omega},\
G(\omega - \tilde{\omega})\rangle\mid\nonumber\\
&\quad=\ \mid\langle u\cdot\nabla\omega - v\cdot\nabla\omega +
v\cdot\nabla\omega - v\cdot\nabla\tilde{\omega},\ G(\omega -
\tilde{\omega})\rangle\mid\nonumber\\
&\quad\leq\ \mid\langle (u-v)\cdot\nabla\omega, \ G(\omega -
\tilde{\omega})\rangle\mid\ + \mid\langle v\cdot\nabla(\omega -
\tilde{\omega}),\ G(\omega -
\tilde{\omega})\rangle\mid.\label{e4.40}
\end{align}
Let us denote $\omega^{*} = \omega - \tilde{\omega}$. Then with the
help of \eqref{e4.10} we obtain
\begin{align}\label{e4.41}
\mid\langle v\cdot\nabla(\omega - \tilde{\omega}),\ G(\omega -
\tilde{\omega})\rangle\mid\ &\leq\ C\parallel v\parallel_{m}\
Q_{p}(\omega - \tilde{\omega}).
\end{align}
Since $\Div(u-v) = 0$, we have
\begin{align}
\mid\langle (u-v)\cdot\nabla\omega, \ G(\omega -
\tilde{\omega})\rangle\mid\ &=\
\mid\langle\nabla\cdot((u-v)\otimes\omega), \
G(\omega^{*})\rangle\mid\nonumber\\
&=\ \mid\langle\nabla\cdot((u-v)\otimes\omega),\
\mid\omega^{*}\mid^{p-2}\omega^{*} \rangle\mid\nonumber.
\end{align}
Integrating by parts we get,
\begin{align}
\mid\langle (u-v)\cdot\nabla\omega, \ G(\omega -
\tilde{\omega})\rangle\mid\ &=\ \mid\langle (u-v)\otimes\omega ,\
\nabla(\mid\omega^{*}\mid^{p-2}\omega^{*}) \rangle\mid\nonumber\\
&\leq\langle\mid(u-v)\otimes\omega\mid, \
\mid\nabla(\mid\omega^{*}\mid^{p-2}\omega^{*})\mid\rangle\nonumber\\
&\leq\langle\mid(u-v)\otimes\omega\mid, \
\mid\omega^{*}\mid^{p-2}\mid\nabla\omega^{*}\mid\rangle.\nonumber
\end{align}
Now using the H\"{o}lder's inequality we obtain,
\begin{align}\label{e4.42}
\mid\langle (u-v)\cdot\nabla\omega, \ G(\omega -
\tilde{\omega})\rangle\mid\ &\leq\parallel
(u-v)\otimes\omega\parallel_{q'} \
\parallel\ \mid\omega^{*}\mid^{p-2}\mid\nabla\omega^{*}\mid\
\parallel_{q},
\end{align}
where $\frac{1}{q} + \frac{1}{q'} = 1$.

Using \eqref{e4.7} and H\"{o}lder's inequality one more time, we
have
\begin{align}\label{e4.43}
\mid\langle (u-v)\cdot\nabla\omega, \ G(\omega -
\tilde{\omega})\rangle\mid\ &\leq\ C\parallel u-v\parallel_{r} \
\parallel\omega\parallel_{m}\
\parallel\omega^{*}\parallel_{r}^{(p-2)/2} \
Q_{p}(\omega^{*})^{1/2},
\end{align}
where $\frac{1}{q'} = \frac{1}{r} + \frac{1}{m}$ and $r =
\frac{mp}{m-2}$.

Notice that if $K$ is the Biot-Savart kernel then $u-v = K*\omega -
K*\tilde{\omega} = K*(\omega - \tilde{\omega}) = K*\omega^{*}$.
Hence using the Lemma 2.7, \eqref{e2.6} we get from \eqref{e4.43}
\begin{align}
&\mid\langle (u-v)\cdot\nabla\omega, \ G(\omega -
\tilde{\omega})\rangle\mid\nonumber\\
&\quad\leq\ C\
(\parallel\omega^{*}\parallel_{\mathrm{L}^{1}}+\parallel\omega^{*}\parallel_{\mathrm{L}^{r}})\
\parallel\omega\parallel_{m} \
\parallel\omega^{*}\parallel_{r}^{(p-2)/2} \
Q_{p}(\omega^{*})^{1/2}\nonumber\\
&\quad=\ C\parallel\omega^{*}\parallel_{r}^{p/2} \
\parallel\omega\parallel_{m} \ Q_{p}(\omega^{*})^{1/2}\nonumber \\
&\quad\quad+ C\parallel\omega^{*}\parallel_{\mathrm{L}^{1}} \
\parallel\omega^{*}\parallel_{r}^{(p-2)/2} \ \parallel\omega\parallel_{m} \
Q_{p}(\omega^{*})^{1/2}.\label{e4.44}
\end{align}
With the help of Lemma 2.6, equation \eqref{e4.44} yields
\begin{align}
&\mid\langle (u-v)\cdot\nabla\omega, \ G(\omega -
\tilde{\omega})\rangle\mid\nonumber\\
&\quad\leq\ C\parallel\omega\parallel_{m} \ Q_{p}(\omega^{*}) + \
C\parallel\omega^{*}\parallel_{\mathrm{L}^{1}} \
\parallel\omega\parallel_{m} \
Q_{p}(\omega^{*})^{(p-1)/p}\nonumber\\
&\quad=\ C\parallel\omega\parallel_{m} \ Q_{p}(\omega
-\tilde{\omega}) + \ C\parallel\omega -
\tilde{\omega}\parallel_{\mathrm{L}^{1}} \
\parallel\omega\parallel_{m} \ Q_{p}(\omega
-\tilde{\omega})^{1/p'}.\label{e4.45}
\end{align}
Thus substituting the results from \eqref{e4.41} and \eqref{e4.45}
in \eqref{e4.40} we have
\begin{align}
\mid\langle u\cdot\nabla\omega - v\cdot\nabla\tilde{\omega},\
G(\omega - \tilde{\omega})\rangle\mid &\leq\ C\ (\ \parallel  v
\parallel_{m} + \parallel\omega\parallel_{m}\ ) \ Q_{p}(\omega
-\tilde{\omega})\nonumber \\ &\quad + \ C\parallel\omega -
\tilde{\omega}\parallel_{\mathrm{L}^{1}} \
\parallel\omega\parallel_{m} \ Q_{p}(\omega
-\tilde{\omega})^{1/p'},\label{e4.46}
\end{align}
where $C$ is a positive constant depending upon $m$ and $p$.

Next we estimate the third term of the equation \eqref{e4.38}. We
notice that
\begin{align}
&\mid\langle\omega\cdot\nabla u - \tilde{\omega}\cdot\nabla v,\
G(\omega - \tilde{\omega})\rangle\mid\nonumber\\
&\quad=\ \mid\langle\omega\cdot\nabla u - \tilde{\omega}\cdot\nabla
u + \tilde{\omega}\cdot\nabla u - \tilde{\omega}\cdot\nabla v, \
G(\omega - \tilde{\omega})\rangle\mid\nonumber\\
&\quad\leq\ \mid\langle (\omega - \tilde{\omega})\cdot u, \ G(\omega
- \tilde{\omega})\rangle\mid + \mid\langle\tilde{\omega}\cdot\nabla
(u-v),\ G(\omega - \tilde{\omega})\rangle\mid.\label{e4.47}
\end{align}
Here we proceed in the similar way as before to get
\begin{align}\label{e4.48}
\mid\langle (\omega - \tilde{\omega})\cdot u, \ G(\omega -
\tilde{\omega})\rangle\mid\ \leq\ C\parallel u\parallel_{m} \
Q_{p}(\omega - \tilde{\omega}),
\end{align}
and
\begin{align}
\mid\langle\tilde{\omega}\cdot\nabla (u-v),\ G(\omega -
\tilde{\omega})\rangle\mid\ &\leq\
C\parallel\tilde{\omega}\parallel_{m} \ Q_{p}(\omega
-\tilde{\omega})\nonumber\\
&\quad+\ C\parallel\omega - \tilde{\omega}\parallel_{\mathrm{L}^{1}}
\
\parallel\tilde{\omega}\parallel_{m} \ Q_{p}(\omega -\tilde{\omega})^{1/p'}.
\end{align}
Thus \eqref{e4.47} yields
\begin{align}
\mid\langle\omega\cdot\nabla u - \tilde{\omega}\cdot\nabla v,\
G(\omega - \tilde{\omega})\rangle\mid &\leq\ C\ (\ \parallel  u
\parallel_{m} + \parallel\tilde{\omega}\parallel_{m}\ ) \ Q_{p}(\omega
-\tilde{\omega})\nonumber\\
&\quad+\ C\parallel\omega - \tilde{\omega}\parallel_{\mathrm{L}^{1}}
\
\parallel\tilde{\omega}\parallel_{m} \ Q_{p}(\omega
-\tilde{\omega})^{1/p'},\label{e4.50}
\end{align}
where $C$ is a positive constant depending upon $m$ and $p$.

Hence \eqref{e4.39}, \eqref{e4.46} and \eqref{e4.50} yield the
desired result \eqref{e4.37} from \eqref{e4.38}.
\end{proof}


\end{document}